\documentclass{article}

\usepackage{amsfonts}
\usepackage{amssymb}
\usepackage{enumerate}
\usepackage{amsmath}
\usepackage[T1]{fontenc}
\usepackage{graphicx}

\newtheorem{ttt}{Theorem}[section]
\newtheorem{llll}[ttt]{Lemma}
\newtheorem{ccc}[ttt]{Claim}
\newtheorem{eee}[ttt]{Example}
\newtheorem{fff}[ttt]{Fact}
\newtheorem{rrr}[ttt]{Remark}
\newtheorem{sss}[ttt]{Statement}
\newtheorem{ddd}[ttt]{Definition}
\newtheorem{qqq}[ttt]{Question}
\newtheorem{cccc}[ttt]{Corollary}
\newtheorem{nnn}[ttt]{Notation}
\newtheorem{ppp}[ttt]{Problem}
\newtheorem{ccccc}[ttt]{Conjecture}

\newcommand{\beq}{\begin{equation} }

\newcommand{\bt}{\begin{ttt}}
\newcommand{\bl}{\begin{llll}}
\newcommand{\bc}{\begin{ccc}}
\newcommand{\bex}{\begin{eee}}
\newcommand{\bfa}{\begin{fff}}
\newcommand{\br}{\begin{rrr}\upshape}
\newcommand{\bst}{\begin{sss}}
\newcommand{\bd}{\begin{ddd}\upshape}
\newcommand{\bq}{\begin{qqq}}
\newcommand{\bnn}{\begin{nnn}}
\newcommand{\bpr}{\begin{ppp}}
\newcommand{\bcor}{\begin{cccc}}
\newcommand{\bcon}{\begin{ccccc}}

\newcommand{\eeq}{\end{equation}}
\newcommand{\et}{\end{ttt}}
\newcommand{\el}{\end{llll}}
\newcommand{\ec}{\end{ccc}}
\newcommand{\eex}{\end{eee}}
\newcommand{\efa}{\end{fff}}
\newcommand{\er}{\end{rrr}}
\newcommand{\est}{\end{sss}}
\newcommand{\ed}{\end{ddd}}
\newcommand{\eq}{\end{qqq}}
\newcommand{\ecor}{\end{cccc}}
\newcommand{\econ}{\end{ccccc}}
\newcommand{\enn}{\end{nnn}}
\newcommand{\epr}{\end{ppp}}

\newcommand{\bp}{\noindent\textbf{Proof. }}
\newcommand{\ep}{\hspace{\stretch{1}}$\square$\medskip}

\newcommand{\lab}[1]{\label{#1}}

\newcommand{\ZZ}{\mathbb{Z}}

\newcommand{\RR}{\mathbb{R}}
\renewcommand{\SS}{\mathbb{S}}

\newcommand{\al}{\alpha}

\newcommand{\iG}{\mathcal{G}}

\newcommand{\sm}{\setminus}

\newcommand{\beeq}{\begin{equation}}
\newcommand{\eeeq}{\end{equation}}
\def\su{\subset}

\newcommand{\dist}{\mathrm{dist}}

\numberwithin{equation}{section}

\title{The structure of continuous rigid functions of two variables}

\author{Rich\'ard Balka\thanks{Partially supported by
    Hungarian Scientific Foundation grant no.~72655.} \\
\\
E\"otv\"os Lor\'and University\\
Department of Analysis\\
P\'az\-m\'any P. s. 1/c, H-1117, Budapest, Hungary\\
balkar@cs.elte.hu \\
\\
and\\
\\
M\'arton Elekes\thanks{Partially supported by
    Hungarian Scientific Foundation grants no.~49786, 61600, and 72655.}\\
\\
Alfr\'ed R\'enyi Institute of Mathematics\\
Hungarian Academy of Sciences\\
P.O. Box 127, H-1364 Budapest, Hungary\\
emarci@renyi.hu\\
www.renyi.hu/ $\tilde{}$ emarci\\
and\\
E\"otv\"os Lor\'and University\\
Department of Analysis\\
P\'az\-m\'any P. s. 1/c, H-1117, Budapest, Hungary\\}

\begin{document}

\maketitle

\insert\footins{\footnotesize{MSC codes: Primary 26A99 Secondary 39B22, 39B52,
    39B72, 51M99}}
\insert\footins{\footnotesize{Key Words: rigid, functional equation,
    transformation, exponential}}

\begin{abstract}
A function $f:\RR^n \to \RR$ is called \emph{vertically rigid} if
$graph(cf)$ is isometric to $graph (f)$ for all $c \neq 0$.
In \cite{BE} we settled Jankovi\'c's conjecture by showing that a continuous
function $f:\RR \to \RR$ is vertically rigid if and only if it is of the form
$a+bx$ or $a+be^{kx}$ ($a,b,k \in \RR$). Now we prove that a continuous
function $f:\RR^2 \to \RR$ is vertically rigid if and only if after a suitable
rotation around the $z$-axis $f(x,y)$ is of the form
$a + bx + dy$, $a + s(y)e^{kx}$ or $a + be^{kx} + dy$ ($a,b,d,k \in \RR$, $k
\neq 0$, $s : \RR \to \RR$ continuous). The problem remains open in higher
dimensions.
\end{abstract}

\section{Introduction}

An easy calculation shows that the exponential function $f(x) = e^x$ has the
somewhat `paradoxical' property that $cf$ is a translate of $f$ for
every $c>0$. It is also easy to see that every function of the form
$a+be^{kx}$ shares this property. Moreover, for every function of the form
$f(x) =
a+bx$ the graph of $cf$ is isometric to the graph of $f$. In \cite{CCR} Cain,
Clark and Rose introduced the notion of vertical rigidity, which we now
formulate for functions of several variables.

\bd
A function $f:\RR^n \to \RR$ is called \emph{vertically rigid}, if $graph(cf)$
is isometric to $graph (f)$ for all $c \in (0,\infty)$. (Clearly, $c \in \RR
\sm \{0\}$ would be the same.)
\ed

Then D.~Jankovi\'c formulated the following conjecture (see \cite{CCR}).
\bcon
\textbf(D.~Jankovi\'c)
A continuous function $f : \RR \to \RR$ is vertically rigid if and only if it
is of the form $a+bx$ or $a+be^{kx}$ ($a,b,k \in \RR$, $k \neq 0$).
\econ

This conjecture, and more, was proved in \cite{BE}.

\bt
\lab{t:Jan}
Jankovi\'c's conjecture holds. (It is actually enough to assume that $f$ is
vertically rigid for an uncountable set $C$, see Definition \ref{d:C} below.)
\et

Later C.~Richter gave generalisations of this theorem in various directions,
see \cite{Ri}.

The main goal of the present paper is to give a complete description of
the continuous vertically rigid functions of two variables.

\bt
\textbf{(Main Theorem)}
A continuous function $f : \RR^2 \to \RR$ is vertically rigid
if and only if after a suitable rotation around the $z$-axis $f(x,y)$
is of the form $a + bx + dy$, $a + s(y) e^{kx}$ or $a + b e^{kx} + dy$
($a,b,d,k \in \RR$, $k \neq 0$, $s : \RR \to \RR$ continuous).
\et

As these classes look somewhat ad hoc, we do not even have conjectures in
higher dimensions.

\bpr
Characterise the continuous vertically rigid functions of $n$ variables for $n
\ge 3$.
\epr

In fact, for the proof of the Main Theorem we need the following technical
generalisations.

\bd
\lab{d:C}
If $C$ is a subset of $(0, \infty)$ and $\iG$ is a set of isometries of
$\RR^3$ then we say that $f$ is vertically rigid \emph{for a set $C \su
(0, \infty)$ via
elements of $\iG$} if for every $c \in C$ there exists a $\varphi \in \iG$
such that $\varphi(graph(f)) = graph (cf)$.

(If we do not mention $C$ or $\iG$ then $C$ is $(0, \infty)$ and $\iG$
is the set of all isometries.)
\ed

\bd
Let us say that a set $C \su (0, \infty)$ \emph{condensates to $\infty$} if
for every $r \in \RR$ the set $C \cap (r ,\infty)$ is uncountable.
\ed

The Main Theorem will immediately follow from the following, in which we just
replace $(0,\infty)$ by a set $C$ condensating to $\infty$.

\bt
\lab{t:mainC}
\textbf{(Main Theorem, technical form)}
Let $C \su (0, \infty)$ be a set condensating to $\infty$. Then a continuous
function $f : \RR^2 \to \RR$ is vertically rigid for $C$ if and only if after
a suitable rotation around the $z$-axis $f(x,y)$ is of the form $a + bx + dy$,
$a + s(y) e^{kx}$ or $a + b e^{kx} + dy$ ($a,b,d,k \in \RR$, $k \neq 0$, $s :
\RR \to \RR$ continuous).
\et

The structure of the proof will be the following.
First we check in Section \ref{s:forms}
that functions of the above forms are rigid. (Of course, they are all
continuous.)
Then we start proving the Main
Theorem in more and more general settings. In Section
\ref{s:translations} first we show that if all the
isometries are horizontal translations then the vertically rigid function
$f(x,y)$ is of the form
$s(y) e^{kx}$ ($k \in \RR$, $k~\neq~0$, $s : \RR~\to~\RR$ continuous).
The punchline here is that we can derive a simple functional equation
from vertical rigidity (some sort of `multiplicativity', see Lemma
\ref{l:add}). Then we
conclude this section by referring to a completely algebraic proof in
\cite{BE} showing that if we allow arbitrary translations then $f(x,y)$ is of
the form $a + s(y) e^{kx}$ ($a, k \in \RR$, $k \neq 0$, $s : \RR \to \RR$
continuous).

Then we start working on the case of general isometries.
The central idea is to consider the set $S_f$  of directions of segments
connecting pairs of points on $graph(f)$ (see Definition \ref{d:S_f}). We
collect the necessary properties of this set in Section \ref{s:propS_f}. The
set $S_f$ has some sort of rigidity in that the transformation $f \mapsto cf$
distorts the shape of it,
but the resulting set has to be isometric to the
original one (see Definition \ref{d:psi} and Remark \ref{r:rigid}). Using
these we determine the possible $S_f$'s in Section \ref{s:possS_f}, then in
Section \ref{s:vertpf} we complete the proof by handling
these cases using various methods.

Finally, in Section \ref{s:open} we collect the open questions.

\section{Functions of these forms are rigid}
\lab{s:forms}

Rotation of the graph around the $z$-axis does not affect vertical rigidity,
so we can assume that $f$ is of the given form without rotations.

Functions of the form $a + bx + dy$ are clearly vertically rigid.

Let now $f(x,y) = a + s(y) e^{kx}$ ($a,k \in \RR$, $k \neq 0$, $s
: \RR \to \RR$ continuous). Then $cf(x,y) = f(x + \frac{\log
c}{k}, y) + a(c-1)$, so $f$ is actually vertically rigid via
translations in the $xz$-plane.

Before checking the third case we need a lemma.

\bl
\label{lem1}
Let $f(x,y)=g(x)+dy$, where $d>0$ and let $c>0$. If
we rotate $graph(f)$ around the $x$-axis by the angle
$\alpha_c=\arctan (cd)-\arctan(d)$ then the intersection of this
rotated graph with the $xy$-plane is the graph of a function of the form $y
= - w_{c,d} g(x) $, where $w_{c,d}>0$ and the map
$c\mapsto w_{c,d}$ is strictly monotone on $(0, \infty)$ for every fixed $d>0$.
\el

\br
By rather easy and short elementary geometric considerations one can check that
for every fixed $d>0$ the map $c \mapsto w_{c,d}$ is positive and real analytic.
It is also very easy
to see geometrically that the limit at $0$ is $\infty$, hence it is not constant,
therefore countable-to-one. This would suffice for all our purposes, but these
arguments are unfortunately very hard to write down rigorously, so we decided
to present a less instructive and longer algebraic proof.
\er

\bp
Using the matrix of the rotation we can write the rotated image
of the point of the graph $(x,y_0,g(x)+dy_0)$ as
\beq \label{e00}
\left(\begin{matrix}
1&0&0\cr
0&\cos \alpha _{c}& -\sin \alpha _{c}\cr
0&\sin \alpha _{c}& \cos \alpha _{c}
\end{matrix}\right)
\left(\begin{matrix}
x\cr
y_0\cr
g(x)+dy_0
\end{matrix}\right)
= \left(\begin{matrix}
x\cr
y_0(\cos \alpha _{c}-d\sin \alpha_{c})-g(x)\sin \alpha _{c}\cr
y_0(\sin \alpha_{c}+d\cos \alpha _{c})+g(x)\cos \alpha _{c}
\end{matrix}\right).
\eeq 
Let us now determine the intersection of the rotated graph
with the $xy$-plane. This right hand side of (\ref{e00}) is in the
$xy$-plane if and only if the third coordinate vanishes, that is,
when $y_0 (\sin \alpha _{c}+d\cos \alpha_{c})+g(x)\cos \alpha
_{c}=0$. This yields 
\beq
\label{e01}
y_{0}=-\frac{\cos \alpha_{c}}{\sin \alpha_{c}+d\cos \alpha _{c}}g(x).
\eeq
In order to
complete the proof of the lemma we have to calculate the
$y$-coordinate of the rotated image of the point
$(x,y_0,g(x)+dy_0)$, which is the second entry of the right hand
side of (\ref{e00}). Hence, using (\ref{e01}), 
\[
y= y_{0} (\cos {\alpha _{c}}-d\sin {\alpha _{c}})-g(x)\sin \alpha
_{c}= -\frac{\cos \alpha _{c}(\cos \alpha _{c}-d\sin \alpha
_{c})}{\sin \alpha _{c}+d\cos \alpha _{c}}g(x)-
\]
\[
g(x)\sin \alpha _{c}= -\frac{\cos ^{2}  \alpha _{c}-d\cos {\alpha
_{c}} \sin {\alpha _{c}}+\sin ^{2} \alpha _{c}+d\cos {\alpha _{c}}
\sin {\alpha _{c}}}{\sin \alpha _{c}+d\cos \alpha _{c}}g(x)=
\]
\[
-\frac{1}{\sin \alpha _{c}+d\cos \alpha _{c}}g(x). \]
Therefore
\[
w_{c,d}=\frac{1}{\sin \alpha _{c}+d\cos \alpha
_{c}}=\left(\sqrt{d^{2}+1}\left(\sin \alpha
_{c}\frac{1}{\sqrt{d^{2}+1}}+\cos \alpha
_{c}\frac{d}{\sqrt{d^{2}+1}}\right)\right)^{-1}.
\]
Using the identity
\beq
\label{e1:43}
\sin \alpha = \frac{\tan
\alpha}{\sqrt{\tan^2 \alpha + 1}} \ \left( \alpha \in (- \pi/2, \pi/2) \right)
\eeq
we obtain $\sin
(\arctan (d))=\frac{d}{\sqrt{d^{2}+1}}$, which easily implies $\cos (\arctan
(d))=\frac{1}{\sqrt{d^{2}+1}}$. (Note that $\arctan(d)\in (-\pi/2,\pi/2)$.) So
\[
w_{c,d}=\left(\sqrt{d^{2}+1}\big(\sin \alpha
_{c}\cos (\arctan (d))+\cos \alpha _{c}\sin (\arctan
(d))\big)\right)^{-1}.
\]
By the formula $\sin (\alpha+\beta)=\sin \alpha \cos
\beta+\cos \alpha \sin \beta$ and the definition of $\alpha _c$
this equals
\[\left(\sqrt{d^{2}+1}\sin  (\alpha _{c}+ \arctan
(d))\right)^{-1}=\left(\sqrt{d^{2}+1}\sin (\arctan
(cd))\right)^{-1}.
\]
Applying \eqref{e1:43} again yields
\[
w_{c,d} = \left(\sqrt{d^{2}+1}\frac{\tan(
\arctan (cd))}{\sqrt{\tan^2 (\arctan(cd)) +
1}}\right)^{-1}=\sqrt{\frac{1}{d^{2}+1}\left(1+\frac{1}{(cd)^{2}}\right)}.
\]
From this form it is easy to see that this function is positive and 
strictly monotone on $(0,\infty)$ for every fixed $d>0$.
\ep

Let now $f(x,y) = a + b e^{kx} + dy$ ($a,b,d,k \in \RR$, $k \neq 0$).
Rescaling the graph in a homothetic way does not affect vertical rigidity, so
we can consider $kf(\frac{x}{k},\frac{y}{k})$ and assume $k=1$.  We may also
assume $b,d \neq 0$, otherwise our function is of one of the previous
forms. Adding a constant, reflecting the graph about the $xz$-plane (needed
only if the signs of $b$ and $d$ differ), multiplying by a nonzero constant, as
well as a translation in the $x$-direction do not affect vertical rigidity,
so by applying these in this order we can assume that $a=0$, $bd>0$, $d=1$,
and $b=1$.

Hence it suffices to check that $f(x,y) = e^x + y$ is vertically
rigid. Let us fix a $c>0$. In every vertical plane of the form
$\{x = x_0\}$ the restriction of $f$ is a straight line of slope
$1$. Rotation around the $x$-axis by angle $\alpha_c=\arctan
(c)-\frac{\pi}{4}$ takes all these lines to lines of slope $c$. By
applying Lemma \ref{lem1} with $g(x)=e^{x}$ and $d=1$, the
intersection of the rotated graph and the $xy$-plane is the graph
of the function $y = - w_{c,1} e^x$.

Now, applying a translation in the $x$-direction we can obtain a
function with still all lines of slope $c$ but now with
intersection with the $xy$-plane of the form $y = -e^x$ (note that
$w_{c,1}>0$). But then we are done, since this function clearly
agrees with $cf$. (The intersection of $graph(f)$ and the
$xy$-plane is of the form $y = -e^x$, and all lines in this graph
are of slope $1$, hence for $graph(cf)$ the intersection is still
$y = -e^x$, and all lines are of slope $c$.) This finishes the
proof of vertical rigidity.

\section{Vertical rigidity via translations}
\lab{s:translations}

\bt
\lab{t:htcont}
Let $C \su (0, \infty)$ be an uncountable set. Then a continuous function $f :
\RR^2 \to \RR$ is
vertically rigid for $C$ \emph{via horizontal translations} if and only if
after a suitable rotation around the $z$-axis $f(x,y)$ is of the form
$s(y) e^{kx}$ ($k \in \RR$, $k \neq 0$, $s : \RR \to \RR$ continuous).
\et

We already checked the easy direction in the previous section.
Before proving the other direction we need some preparation. We
will need the following result, which is Theorem 2.5 in \cite{BE}.

\bt \label{old} Let $f:\RR\rightarrow \RR$ be a continuous
vertically rigid function for an uncountable set $C\subset
(0,\infty)$ via horizontal translations. Then $f$ is of the form
$se^{kx}$ $(s\in \RR, k\in \RR\setminus \{0\}).$ \et

The following lemma will be useful throughout the paper. Sometimes
we will use it tacitly. The easy proof is left to the reader.

\bl
\lab{l:C}
Let $f : \RR^2 \to \RR$ be vertically rigid for $c_0$ via $\varphi_0$ and for
$c$ via $\varphi$. Then $c_0 f$ is vertically rigid for $\frac{c}{c_0}$ via
$\varphi\circ \varphi_0^{-1}$.
\el

From now on we will often use the notation $\vec{x}$ for two-dimensional
(and sometimes three-dimensional) vectors.

\bd
For a function $f:\RR^2\to\RR$ and a set $C \su (0, \infty)$
let $T_{f,C} \su \RR^2$ be the
additive group generated by the set $T' = \{ \vec{t} \in \RR^2 : \exists c \in C \
\forall \vec{x} \in \RR^2 \  f(\vec{x}+\vec{t})=cf(\vec{x}) \}$.
(We will usually simply write $T$ for $T_{f,C}$.)
\ed

\bl
\lab{l:add}
Let $f:\RR^2 \rightarrow \RR$ be a vertically rigid function for a set
$C\subset (0,\infty)$ via horizontal translations such that
$f(\vec{0})=1$. Then
\[
f(\vec{x}+\vec{t}) = f(\vec{x}) f(\vec{t}) \ \ \forall \vec{x} \in \RR^2 \
\forall \vec{t} \in T.
\]
Moreover, $f(\vec{t})>0$ for every $\vec{t} \in T$, and $T'$ is uncountable if
so is $C$.
\el

\bp
By assumption, for every $c \in C$ there exists $\vec{t_c} \in \RR^2$
such that $cf(\vec{x}) = f(\vec{x} + \vec{t_c})$ for every $\vec{x} \in
\RR^2$. Then $\vec{t_c} \in T'$ for every $c \in C$.


Since $T$ is the group generated by $T'$, every $\vec{t} \in T$
can be written as $\vec{t} = \sum_{i=1}^m n_i \vec{t_i}$
($\vec{t_i} \in T', n_i \in \ZZ, i=1, \dots, m$) where $f(\vec{x}
+ \vec{t_i})=c_i f(\vec{x}) \ (\vec{x} \in \RR^2, \ i=1, \dots,
m$).

From these we easily get
\beeq
\lab{2}
f(\vec{x} + \vec{t}) = c_{\vec{t}} f(\vec{x}),
\textrm{ where } c_{\vec{t}} = \prod_{i=1}^m c_i^{n_i}, \ \vec{x} \in
\RR^2, \ \vec{t} \in T.
\eeq
Note that $c_{\vec{t}}>0$
(and also that it is not necessarily a member of $C$). It suffices
to show that $c_{\vec{t}} = f(\vec{t})$ for every $\vec{t} \in T$, but this
follows if we substitute $\vec{x} = \vec{0}$ into (\ref{2}).

Since $f$ is not
identically zero, $\vec{t}_c \neq \vec{t}_{c'}$ whenever $c,c' \in C$
are distinct. Hence $\{ \vec{t}_c : c \in C \}$ is uncountable, so $T'$
is uncountable if so is $C$.
\ep

\bp (Thm. \ref{t:htcont}) If $f$ is identically zero then we are
done, so let us assume that this is not the case.  The class of
continuous vertically rigid functions for some set condensating to
$\infty$ via
horizontal translations, as well as the class of functions of the
form $s(y)e^{kx}$ ($k\in \RR$, $k \neq 0$, $s : \RR \to \RR$ continuous)
are both closed
under horizontal translations and under multiplication by nonzero
constants (by Lemma \ref{l:C}).
Hence we may assume that $f(\vec{0})=1$. Then the previous lemma
yields that $f(\vec{t_1} + \vec{t_2}) = f(\vec{t_1})f(\vec{t_2})$
$(\vec{t_1}, \vec{t_2} \in T)$, and also that $f|_T>0$. Then $g(\vec{t}) =
\log f(\vec{t}) $ is defined for every $\vec{t} \in T$, and $g$ is clearly
additive on $T$.

Let us now consider $\bar{T}$, the closure of $T$, which is clearly an
uncountable closed subgroup of $\RR^2$. It is well-known that every closed
subgroup of $\RR^2$ is a nondegenerate linear image of a group of the form
$G_1 \times G_2$, where $G_1, G_2 \in \{ \{0\}, \ZZ, \RR \}$.
Hence after a suitable rotation around the origin
$\bar{T}$ is either $\RR^2$ or $\RR \times \{0\}$ or $\RR \times r\ZZ$ for
some $r>0$.

\noindent\textbf{Case 1.} $\bar{T} = \RR^2$.

In this case $T \su \RR^2$ is dense. It is well-known that a continuous
additive function on a dense subgroup is of the form $g(x,y) = \alpha x +
\beta y$, $((x,y) \in T)$ for some $\alpha, \beta \in \RR$. But then $f(x,y) =
e^{\alpha x + \beta y}$ on $T$, and by continuity this holds on the whole
plane as well. As the constant $1$ function is not vertically rigid via
horizontal translations, $\alpha = \beta = 0$ cannot hold. By applying a
rotation of angle $\frac{\pi}{2}$ if necessary we may assume that $\alpha \neq
0$. But then by choosing $k = \alpha$, $s(y) = e^{\beta y}$ we are
done.

\noindent\textbf{Case 2.} $\bar{T} = \RR \times \{0\}$.

In this case every $\vec{t}_c$ is of the form $(t_c,0)$, where $t_c \neq 0$ if
$c \neq 1$. (We may assume $1 \notin C$.)

Applying Theorem \ref{old} for every fixed $y$ we obtain that
$f(x,y) = s(y) e^{k_y x}$ ($s(y),k_y \in \RR, k_y \neq 0$). As
$s(y) = f(0,y)$, we get that $s$ is continuous. If $s(y) \neq 0$
then it is not hard to see that $k_y = \frac{\log c}{t_c}$, which
is independent of $y$, so for these $y$'s $k_y = k$ is constant.
But if $s(y) = 0$ then the value of $k_y$ is irrelevant, so it can
be chosen to be the same constant $k$. Hence without loss of
generality $k_y = k$ is constant, and we are done with this case.

\noindent\textbf{Case 3.} $\bar{T} = \RR \times r\ZZ$.

As $T'$ is uncountable, there is an $n \in \ZZ$
so that $T' \cap (\RR \times \{rn\})$
is uncountable. Fix an element $t_{c_0}$ of this set. Then Lemma \ref{l:C}
yields that $c_0 f$ is vertically rigid for an uncountable set via
translations of the form $(t,0)$. Restricting ourselves to these isometries
and $c$'s we are done using Case 2, since every uncountable set in $\RR$
generates a dense subgroup.
\ep

Now we handle the case of arbitrary translations.

\bt
\lab{t:trans}
Let $f:\RR^2 \rightarrow \RR$ be an arbitrary
function that is vertically rigid for a set $C \su (0,\infty)$ via
translations. Then there exists $a\in \RR$ such that $f-a$ is
vertically rigid for the same set $C$ via horizontal translations.
\et

\bp
The obvious modification of \cite[Thm. 2.4]{BE} works, just replace all $x$'s
and $u$'s by vectors.
\ep

This readily implies the following.

\bcor
\lab{c:trans}
Let $C \su (0, \infty)$ be an uncountable set. Then a continuous function $f :
\RR^2 \to \RR$ is
vertically rigid for $C$ \emph{via translations} if and only if
after a suitable rotation around the $z$-axis $f(x,y)$ is of the form
$a+s(y) e^{kx}$ ($a, k \in \RR$, $k \neq 0$, $s : \RR \to \RR$ continuous).
\ecor

\section{The set $S_f$}
\lab{s:propS_f}

Now we start working on the case of arbitrary isometries.

Let $\SS^2 \su \RR^3$ denote the unit sphere.
For a function $f : \RR^2 \to \RR$ let $S_f$ be the set of directions between
pairs of points on the graph of $f$, that is,

\bd
\lab{d:S_f}
\[
S_f = \left\{ \frac{p-q}{|p-q|} \in \SS^2 : p,q \in graph(f),\  p \neq q
\right\}.
\]
\ed

Recall that a \emph{great circle} is a circle line in $\RR^3$ of radius $1$
centered at the origin. We call it \emph{vertical} if it passes through the
points $(0,0,\pm 1)$.

\bl
\lab{l:S}
Let $f : \RR^2 \to \RR$ be continuous. Then
\begin{enumerate}
\item
\lab{symm}
$- S_f = S_f$ (symmetric about the origin),
\item
\lab{poles}
$(0,0,\pm 1) \notin S_f$,
\item
\lab{conn}
$S_f$ is connected,
\item
\lab{arcs}
every great circle containing $(0,0,\pm 1)$ intersects $S_f$ in two
(symmetric) nonempty arcs,
\item
\lab{comp}
$\SS^2 \sm S_f$ has exactly two connected components, one containing $(0,0,1)$
and one containing $(0,0,-1)$.
\end{enumerate}
\el

\bp
(\ref{symm}.) Obvious.

(\ref{poles}.) Obvious, since $f$ is a function.

(\ref{conn}.) $graph(f)$ is homeomorphic to $\RR^2$, hence the squared of it
minus the (2-dimensional) diagonal is a connected set. Since $S_f$ is the
continuous image of this connected set, it is itself connected.

(\ref{arcs}.) The intersection of $S_f$ with such a great circle corresponds
to restricting our attention to distinct pairs of points $(\vec{x_1}, \vec{x_2})
\in \RR^2 \times \RR^2$ so that the segment $[\vec{x_1}, \vec{x_2}]$ is
parallel to a fixed line $L \su \RR^2$. Now, given two such nondegenerate
segments it is easy to move one of them continuously to the other so that
along the way it remains nondegenerate and parallel to $L$. This shows that in
both halves of the great circle (separated by $(0,0,\pm 1)$) $S_f$ is pathwise
connected, hence it is an arc.

(\ref{comp}.) By (\ref{arcs}.) every point of $\SS^2 \sm S_f$ can be connected
with an arc of a vertical great circle either to $(0,0,1)$ or to $(0,0,-1)$ in
$\SS^2 \sm S_f$, hence there are at most two connected components.

Now we show that $(0,0,1)$ and $(0,0,-1)$ are in different ones. It suffices to
show that there exists a Jordan curve in $S_f$ so that $(0,0,1)$ and $(0,0,-1)$
are in the two distinct components of its complement.
Let $\SS^1$ denote the unit circle in $\RR^2 =
\{(x,y,z) : z=0\}$ and let $\gamma :\SS^1 \rightarrow S_{f}$ be given by 
\[
\gamma(\vec{x}) =
\frac{(\vec{x},f(\vec{x}))-(-\vec{x},f(-\vec{x}))}{|(\vec{x},f(\vec{x}))-(-\vec{x},f(-\vec{x}))|}.
\]

In this paragraph the word `component' will refer to the components of $\SS^2
\sm \gamma(\SS^1)$.
One can easily check that $\gamma$ is continuous and injective, hence a Jordan curve.
Moreover, it is clearly in $S_f$, and its intersection with every vertical great
circle is a symmetric pair of points. Therefore every point of $\SS^2 \sm
\gamma(\SS^1)$ can be connected with an arc of a vertical great circle either to
$(0,0,1)$ or to $(0,0,-1)$ in $\SS^2 \sm \gamma(\SS^1)$, hence the union of the
components of $(0,0,1)$ and $(0,0,-1)$ cover $\SS^2 \sm \gamma(\SS^1)$. So 
$(0,0,1)$ and $(0,0,-1)$ are in different components, otherwise
$\SS^2 \sm \gamma(\SS^1)$ would be connected, but the complement of a Jordan
curve in $\SS^2$ has two components.
\ep

The above lemma shows that $S_f$ is something like a `strip around the
sphere'. Now we make this somewhat more precise by defining the top and the
bottom `boundaries' of this strip.

\bd
\lab{d:h}
Let $h : \SS^1 \to \SS^2$ be defined as follows. Every $\vec{x} \in \SS^1$ is
in a unique half great circle connecting $(0,0,1)$ and $(0,0,-1)$. The
intersection of $S_f$ with this great circle is an arc, define $h(\vec{x})$ as
the top endpoint of this arc.
\ed

Clearly, the bottom endpoint of this arc is $-h(-\vec{x})$, so the `top
function bounding the strip $S_f$ is $h(\vec{x})$ and the bottom function is
$-h(-\vec{x})$'. The coordinate functions of $h$ are denoted by
$(h_1,h_2,h_3)$, where $h_3 : \SS^1 \to [-1, 1]$ encodes all information about
$h$.

\bl
\lab{l:h}
Let $f : \RR^2 \to \RR$ be continuous, and $h$ be defined as above. Then
\begin{enumerate}
\item
\lab{notpole}
$h(\vec{x}) \neq (0,0,-1)$ for every $\vec{x} \in \SS^1$
\item
\lab{semicont}
$h$ is lower semicontinuous (in the obvious sense, or equivalently, $h_3$ is
lower semicontinuous)
\item
\lab{convex}
$h$ is convex with respect to
great circles, that is, if $h(\vec{x})$ and $h(\vec{y})$
determine a unique nonvertical great circle (i.e.~there is a subarc of
$\SS^1$ of length~$< \pi$ connecting $\vec{x}$ and $\vec{y}$, and $h(\vec{x}),
h(\vec{y}) \neq (0,0,1)$) then on this subarc $graph(h)$ is bounded from above
by the great circle.
\end{enumerate}
\el

\bp
(\ref{notpole}.)
Obvious by Lemma \ref{l:S} (\ref{poles}.) and (\ref{arcs}.).

(\ref{semicont}.)
We have to check that if $h_3(\vec{x}) > u$ then the same holds in a
neighbourhood of $\vec{x}$. (Note that essentially $h_3$ is defined as a
supremum.) Hence $h_3(\vec{x}) > u$ if and only if there exists a segment
$[\vec{a},\vec{b}] \su \RR^2$ parallel to $\vec{x}$ over which the slope of
$f$ is bigger than $u$. But then by the continuity of $f$ the same holds for
segments close enough to $[\vec{a},\vec{b}]$, in particular to slightly
rotated copies, and we are done.

(\ref{convex}.) It is easy to see that for every $\vec{v}\in
\SS^1$ the slope of $f$ over a segment parallel to $\vec{v}$ is at most
the slope of the vector $h(\vec{v})$. Let $\vec{z} \in
\SS^1$ be an element of the shorter arc connecting $\vec{x}$ and
$\vec{y}$ in $\SS^1$, let $[\vec{a},\vec{b}] \su \RR^2$ be a
segment parallel to $\vec{z}$, and let $P:\RR^{2}\rightarrow \RR$
be the linear map whose graph passes through the origin,
$h(\vec{x})$ and $h(\vec{y})$. (Then $graph(P)$ contains the great
circle determined by $h(\vec{x})$ and $h(\vec{y})$. Moreover, 
the slope of $P$ over any vector parallel to $\vec{x}$ is the slope of
$h(\vec{x})$, and similarly for $\vec{y}$.) We have to
show that the slope of $f$ between $\vec{a}$ and $\vec{b}$ is at most
that of $P$, that is, $f(\vec{b})-f(\vec{a})\leq
P(\vec{b})-P(\vec{a})$. Write $\vec{b} - \vec{a} = \al \vec{x} +
\beta \vec{y}$ for some $\al,\beta > 0$. Then by using the
definition of $P$ and our first observation for the segments
$[\vec{a},\vec{a}+\al \vec{x}]$ and $[\vec{a}+\al
\vec{x},\vec{a}+\al \vec{x}+\beta \vec{y}]$, which are parallel to
$\vec{x}$ and $\vec{y}$, respectively, we get
\[
f(\vec{b})-f(\vec{a})=f(\vec{a}+\al
\vec{x}+\beta \vec{y})-f(\vec{a})=$$ $$ \left(f(\vec{a}+\al
\vec{x}+\beta \vec{y}\right)-f(\vec{a}+\al
\vec{x}))+\left(f(\vec{a}+\al \vec{x})-f(\vec{a})\right)\leq 
\]
\[
\left(P(\vec{a}+\al \vec{x}+\beta \vec{y})-P(\vec{a}+\al
\vec{x})\right)+\left(P(\vec{a}+\al
\vec{x})-P(\vec{a})\right)=P(\vec{b})-P(\vec{a}).
\]
\ep

\section{Determining the possible $S_f$'s}
\lab{s:possS_f}

\bd
\lab{d:psi}
For $c > 0$ let $\psi_c: \SS^2 \to \SS^2$ denote the map that `deforms $S_f$
according to the map $c \mapsto cf$', that is,
\[
\psi_c((x,y,z)) = \frac{(x,y,cz)}{|(x,y,cz)|} \ \ ((x,y,z) \in \SS^2).
\]
\ed

\br \lab{r:rigid} Let $\varphi_c$ be the isometry mapping
$graph(f)$ onto $graph(cf)$. Every isometry $\varphi$ is of the
form $\varphi^{trans} \circ \varphi^{ort}$, where $\varphi^{ort}$
is an orthogonal transformation and $\varphi^{trans}$ is a
translation. Moreover, if $\varphi$ is orientation-preserving then
$\varphi^{ort}$ is a rotation around a line passing through the
origin.  A key observation is the following: The vertical rigidity
of $f$ for $C$ implies that $\psi_c(S_f) =\varphi_c^{ort}(S_f)$
for every $c \in C$. \er

Now we prove the main theorem of this section. For the definition of $h_3$ see
the previous section.

\bt
\lab{t:cases}
Let $C \su (0, \infty)$ be a set condensating to $\infty$, and let $f : \RR^2
\to \RR$ be a continuous function vertically rigid for $C$.
Then one of the following holds.
\begin{itemize}

\item
\textbf{Case A.} There is a vertical great circle that intersects $S_f$ in only
two points.

\item
\textbf{Case B.} $S_f = \SS^2 \sm \{ (0,0,1), (0,0,-1) \}$.

\item
\textbf{Case C.}
There exists an $\vec{x_0} \in \SS^1$ such that $h_3(\vec{x_0}) = 0$ and
$h_3(\vec{x}) = 1$ for every $\vec{x} \neq \vec{x_0}$, that is, $S_f$ is
`$\SS^2$ minus two quarters of a great circle'.

\item
\textbf{Case D.}
There exists a closed interval $I$ in $\SS^1$ with $0 < length(I) < \pi$ such
that $h_3(\vec{x}) = 0$ if $\vec{x} \in I$, and
$h_3(\vec{x}) = 1$ if $\vec{x} \notin I$, that is, $S_f$ is
`$\SS^2$ minus two spherical triangles'.

\end{itemize}
\et

\bp In this proof the word `component' will refer to the
components of $\SS^2\setminus S_{f}$. We separate two cases
according to whether $h_3 \ge 0$ everywhere or not.

First let us suppose that there exists a
$\vec{x} \in \SS^1$ such that $h_3(\vec{x}) < 0$. This implies that there is a
vertical great circle containing two arcs, one in the top component
connecting $(0,0,1)$ with $\SS^1$ and even crossing it, and an other
one (the symmetric pair in the bottom component)
running from the `South Pole to the Equator' and even
above. But then considering geometrically the action of $\psi_c$ one can
easily check that if we choose larger and larger $c$'s (tending to $\infty$)
then we obtain that $\psi_c(S_f)$ contains in the two components two
symmetrical arcs on the same great circle which are only leaving out two small
gaps of length tending to $0$. But then by Remark \ref{r:rigid} $S_f$ also
contains two such arcs in the two components on some (not necessarily
vertical) great circle, hence the distance of the components is $0$.

Let $\vec{p_n}$ and $\vec{q_n}$ be sequences in the top and bottom component,
respectively, so that $\dist(\vec{p_n}, \vec{q_n}) \to 0$. By compactness we may assume
$\vec{p_n}, \vec{q_n} \to \vec{p} \in \SS^2$. We claim that $\vec{p_n} \to \vec{p}$ implies $\vec{p} \neq
(0,0,-1)$. (And similarly $\vec{q_n} \to \vec{p}$ implies $\vec{p} \neq (0,0,1)$.) Indeed, let
$\vec{x}_n \in \SS^1$ be so that $\vec{x}_n$ and $\vec{p_n}$ lay on the same
vertical great circle, and similarly, let $\vec{x} \in \SS^1$ and $\vec{p}$ lay on
the same vertical great circle. Then $\vec{x}_n \to \vec{x}$, and using the
fact $h(\vec{x}) \neq (0,0,-1)$ and the lower semicontinuity of $h$ at
$\vec{x}$ (Lemma \ref{l:h} (\ref{notpole}.) and (\ref{semicont}.)) we are done.

Using the lower semicontinuity of $h$ at $\vec{x}$ again (and $\vec{p_n} \to \vec{p}$)
we get that
$h(\vec{x})$ cannot be above $\vec{p}$. Similarly, $-h(-\vec{x})$ cannot be below
$\vec{p}$. But $h(\vec{x})$ is always above $-h(-\vec{x})$, so the only option is
$h(\vec{x}) = -h(-\vec{x})$, hence there is a vertical
great circle whose intersection with $S_f$ is just a (symmetric) pair of
points, so Case A holds, and hence we are done with the first half of the
proof.

Now let us assume that $h_3 \ge 0$ everywhere. First we prove that
$h_3(\vec{x}) \in \{0,1\}$ for Lebesgue almost every $\vec{x} \in
\SS^1$. Indeed, fix an arbitrary $c \in C \sm \{1\}$.  By rigidity the (equal)
measure of the
two components remains the same after applying $\psi_c$. Since $h_3 \ge 0$, the
intersection of the top component with the vertical great circle containing an
$\vec{x}$ shrinks if $c>1$ and grows if $c<1$, unless $h_3(\vec{x}) = 0$ or
$1$. Hence we are done, since the measure of the top component can be
calculated from the lengths of these arcs.

Now we show that $\{\vec{x} : h_3(\vec{x}) = 0 \}$ is either empty, or a pair
of points of the form $\{\vec{x_0}, -\vec{x_0}\}$, or a closed interval in
$\SS^1$ (possibly degenerate or the whole $\SS^1$). So we have to show that if
$\vec{x}, \vec{y} \in \SS^1$ are so that the shorter
arc connecting them is shorter than $\pi$, and $h_3(\vec{x}) = h_3(\vec{y}) =
0$ then $h_3(\vec{z}) = 0$ for every $\vec{z}$ in this arc. But $h_3(\vec{z})
\ge 0$ by assumption, and $h_3(\vec{z})\le 0$ by the convexity of $h$ applied
to $h(\vec{x}) = \vec{x}$ and $h(\vec{y}) = \vec{y}$.
The fact that the endpoints are also contained in $\{\vec{x} : h_3(\vec{x}) =
0 \}$ easily follows from the semicontinuity.

If $\{\vec{x} : h_3(\vec{x}) = 0 \}$ is a symmetrical pair
of points or a closed interval of length at least $\pi$ then it is easy to see
that Case A holds. Hence we may assume that it is empty, or a singleton, or a
closed interval $I$ with $0 < length(I) < \pi$.

\noindent\textbf{Case 1.} $\{\vec{x} : h_3(\vec{x}) = 0 \} = \emptyset$.

In this case, $h_3 > 0$ everywhere, and hence $h_3 = 1$ almost
everywhere. Therefore one can easily see (using the convexity) that $h_3 = 1$
everywhere but possibly at at most two points of the form $\{\vec{x_0},
-\vec{x_0}\}$. We claim that actually $h_3 = 1$ everywhere. We know already
that $S_f$ is $\SS^2$ minus two symmetric arcs on the same vertical great
circle. The arcs contain $(0,0,1)$ and $(0,0,-1)$, respectively, and they do
not reach the `Equator', since $h_3 > 0$. Let us fix an arbitrary $c \in C \sm
\{ 1 \}$.  By rigidity the
(equal) length of the arcs should not change when applying $\psi_c$,
but it clearly changes, a contradiction.

Hence $S_f = \SS^2 \sm\{(0,0,1), (0,0,-1)\}$, so Case B holds.

\noindent\textbf{Case 2.} $\{\vec{x} : h_3(\vec{x}) = 0 \}$ is a singleton.

Let $\{\vec{x_0}\} = \{\vec{x} : h_3(\vec{x}) = 0 \}$. Similarly as above,
$h_3 = 1$ almost everywhere. Then convexity easily implies that $h_3(\vec{x})
= 1$ whenever $\vec{x} \notin \{\vec{x_0}, -\vec{x_0}\}$.
Again similarly, the length of the
arcs is unchanged by $\psi_c$ only if $h_3(\vec{-x_0}) = 1$, so $S_f$ is
$\SS^2$ minus two symmetric quarter arcs starting from the `Poles' on a
vertical great circle, so Case C holds.

\noindent\textbf{Case 3.} $\{\vec{x} : h_3(\vec{x}) = 0 \}$ is a closed interval
in $\SS^1$ with $0 < length(I) < \pi$.

Let $I = \{\vec{x} : h_3(\vec{x}) = 0 \}$.
As $h_3 = 0$ or $1$ almost everywhere, convexity readily implies that $h_3 =
1$ on $\SS^1 \sm I$. Hence $S_f$ is `$\SS^2$ minus two spherical
triangles', and Case D holds.

This concludes the proof.
\ep

\section{The end of the proof}
\lab{s:vertpf}

Now we complete the proof of the technical form of the
Main Theorem. We repeat the statement here.

\bt
(Main Theorem, technical form)
Let $C \su (0, \infty)$ be a set condensating to $\infty$. Then a continuous
function $f : \RR^2 \to \RR$ is vertically rigid for $C$ if and only if after
a suitable rotation around the $z$-axis $f(x,y)$ is of the form $a + bx + dy$,
$a + s(y) e^{kx}$ or $a + b e^{kx} + dy$ ($a,b,d,k \in \RR$, $k \neq 0$, $s :
\RR \to \RR$ continuous).
\et

\bp
By Theorem \ref{t:cases} it suffices to consider Cases A-D.

\noindent\textbf{Case A.} There is a vertical great circle that intersects
$S_f$ in only two points.

We may assume using a suitable rotation around the $z$-axis that the vertical
great circle is in the $yz$-plane, hence $f(x,y)$ is of the form $g(x) + dy$.
The continuity of $f$ implies that $g$ is also continuous. 

\textbf{Subcase A1.} $d = 0$.

Let $c \in C$ be fixed, and let $\varphi_c$ be the corresponding isometry. The
graph of $cf$ is invariant under translations parallel to the $y$-axis. As the
same holds for $f$, by rigidity, $cf$ is also invariant under translations
parallel to the $\varphi_c$-image of the $y$-axis. If these two directions are
nonparallel, then $graph(cf)$ is a plane, and hence so is $graph(f)$, so we
are done since $f(x,y)$ is of the form $a + bx$ (note that there is no `$+
dy$' since $f$ does not depend on $y$). Therefore we may assume that all lines
parallel to the $y$-axis are taken to lines parallel to the $y$-axis, but then
all planes parallel to the $xz$-plane are taken to planes parallel to the
$xz$-plane. But this shows (by considering the intersections of the graphs
with the $xz$-plane) that $g$ is vertically
rigid for $c$, hence by Theorem \ref{t:Jan}
$g(x)$ is of the from $a + bx$ or $a + b e^{kx}$ ($a,b,k \in \RR$, $k \neq
0$), and we are done.

\textbf{Subcase A2.} $d \neq 0$.

We may assume that $d>0$, since otherwise we may consider $-f$.

For every $c \in C$ let $\varphi_c$ be the corresponding isometry.
We claim that we may assume that all these are
orientation-preserving. If $\{c \in C : \varphi_c \textrm{ is
orientation-preserving} \}$ condensates to $\infty$ then we are
done by shrinking $C$, otherwise we may assume that they are all
orientation-reversing (note that if we split $C$ into two pieces
then at least one of them still condensates to $\infty$). Let us
fix a $c_0 \in C$ and consider $c_0f$ instead of $f$. By Lemma
\ref{l:C} this function is rigid for an uncountable set with all
isometries orientation-preserving, and if it is of the desired
form then so is $f$, so we are done.

We may assume $1 \notin C$. Let us fix a $c \in C$. Similarly as
in the previous subcase, we may assume that lines parallel to
$(0,1,d)$ are taken to lines parallel to $(0,1,cd)$ as follows.
The special form of $f$ implies that $graph(f)$ is invariant under
translations in the $(0,1,d)$-direction, hence $graph(cf)$ is
invariant under translations in the $(0,1,cd)$-direction,
moreover, by rigidity, $graph(cf)$ is also invariant under
translations parallel to the $\varphi_c$-image of the lines of
direction $(0,1,d)$. If these two latter directions do not
coincide then $graph(cf)$ is a plane, and we are done.

Therefore the image of every line parallel to $(0,1,d)$ is a line
parallel to $(0,1,cd)$ under the orientation-preserving isometry
$\varphi_c$. As in Remark \ref{r:rigid}, write $\varphi_c = \varphi_c^{trans}
\circ \varphi_c^{ort}$, where $\varphi_c^{ort}$ is a rotation about
a line containing the origin and $\varphi_c^{trans}$ is a translation.
Since the translation does not affect directions, the rotation
$\varphi_c^{ort}$ takes the direction $(0,1,d)$ 
to the nonparallel direction $(0,1,cd)$ ($d \neq
0$), therefore the axis of the rotation has to be orthogonal to the plane
spanned by these two directions. Hence the axis has to be the
$x$-axis. Moreover, the angle of the rotation is easily seen to be 
$\arctan(cd) - \arctan(d)$.

We now show that we may assume that $\varphi_c^{trans}$ is a horizontal translation. Decompose the translation as $\varphi_c^{trans} = \varphi_c^{\vec{u}} \circ \varphi_c^{\vec{v}}$, 
where $\varphi_c^{\vec{v}}$ is a horizontal translation and $\varphi_c^{\vec{u}}$ is a translation in the $(0,1,cd)$-direction. Since $\varphi_c^{ort} (graph(f))$ is invariant under translations in the $(0,1,cd)$-direction, so is $\varphi_c^{\vec{v}} \circ \varphi_c^{ort} (graph(f))$, hence
\[
\varphi_c^{\vec{v}} \circ \varphi_c^{ort} (graph(f)) = \varphi_c^{\vec{u}} \circ \varphi_c^{\vec{v}} \circ \varphi_c^{ort} (graph(f)) = \varphi_c (graph(f)) = graph(cf),
\]
so we can assume $\varphi_c = \varphi_c^{\vec{v}} \circ \varphi_c^{ort}$, and we are done.

We will now complete the proof of this subcase by showing that the function $-\frac{1}{d} g$ is rigid for an uncountable set. Indeed, this suffices by Theorem \ref{t:Jan} and by the special form of $f$.

Let us denote the $xy$-plane by $\{z=0\}$ and consider the intersection of both sides of the equation $\varphi_c (graph(f)) = graph(cf)$ with $\{z=0\}$. 
On the one hand, $\{z=0\} \cap \varphi_c (graph(f)) = 
\{z=0\} \cap \varphi_c^{\vec{v}} \circ \varphi_c^{ort} (graph(f)) =
\varphi_c^{\vec{v}} ( \{z=0\} \cap \varphi_c^{ort} (graph(f))) =
\varphi_c^{\vec{v}} \left( graph \left( -w_{c,d} g \right)\right) =
\varphi_c^{\vec{v}} \left( graph \left( (w_{c,d}d)(-\frac{1}{d} g \right)\right)$,
where we used the fact that $\varphi_c^{\vec{v}}$ is horizontal and Lemma \ref{lem1}.
On the other hand, it is easy to see that $\{z=0\} \cap graph(cf) = graph\left( -\frac{1}{d} g \right)$. Therefore $graph\left( -\frac{1}{d} g \right) = \varphi_c^{\vec{v}} \left( graph \left( (w_{c,d}d)(-\frac{1}{d} g \right)\right)$ and hence 
$-\frac{1}{d} g$
is rigid for $w_{c,d}d$ for every $c>0$. The map $c \mapsto
w_{c,d}d$ is strictly monotone for every fixed $d$, hence the range of $C$ is
uncountable. So $-\frac{1}{d} g$ is rigid for an uncountable set,
and we are done.

\noindent\textbf{Case B.} $S_f = \SS^2 \sm\{(0,0,1), (0,0,-1)\}$.

So $S_f$ is invariant under every $\psi_c$, and hence so is under every
$\varphi_c^{ort}$. Then clearly $\varphi_c^{ort}((0,0,1)) = (0,0,1)$ or
$\varphi_c^{ort}((0,0,1)) = (0,0,-1)$ for every $c \in C$. By the same
argument as above we can assume that the former holds for
every $c \in C$. Using the argument again we can assume
that all $\varphi_c$'s are orientation-preserving. But then each of these is a
rotation around the $z$-axis followed by a translation, in other words, an
orientation-preserving transformation in the $xy$-plane followed by a
translation in the $z$-direction. An orientation-preserving transformation in
the plane is either a translation or a rotation. If it is a translation for
every $c$ then we are done by Corollary \ref{c:trans}. So let us assume that
there exists a $c$ such that $\varphi_c$ is a proper rotation around $\vec{x}
\in \RR^2$ followed by a vertical translation. We claim that then $f$ is
constant, which will contradict that $S_f$ is nearly the full sphere, finishing
the proof of this case. We will actually show that
$f$ is constant on every closed disc $B(\vec{x},R)$ centered at
$\vec{x}$. Indeed, consider $\max_{B(\vec{x},R)} f - \min_{B(\vec{x},R)}
f$. This is unchanged by the rotation around $\vec{x}$ as well as by the
vertical translation, hence by $\varphi_c$. But the map $f \mapsto cf$
multiplies this amount by $c \neq 1$, so the only option is
$\max_{B(\vec{x},R)} f - \min_{B(\vec{x},R)} f = 0$, and we are done.

\noindent\textbf{Case C.} There exists an $\vec{x_0} \in \SS^1$ such that
$h_3(\vec{x_0}) = 0$ and $h_3(\vec{x}) = 1$ for every $\vec{x} \neq
\vec{x_0}$, that is, $S_f$ is `$\SS^2$ minus two quarters of a great circle'.

So $S_f$ is invariant under every $\psi_c$, and hence so is under every
$\varphi_c^{ort}$. Hence $\varphi_c^{ort}$ maps $(0,0,1)$ to one of the four
endpoints of the two arcs. Therefore we can assume by splitting $C$ into four
pieces according to the image of $(0,0,1)$ and applying Lemma \ref{l:C}
that $(0,0,1)$ is a fixed point of every $\varphi_c^{ort}$. But then the two
arcs are also fixed, and actually $\varphi_c^{ort}$ is the identity. Hence
every $\varphi_c$ is a translation, and we are done by Corollary
\ref{c:trans}.

\noindent\textbf{Case D.} There exists a closed interval $I$ in $\SS^1$ with
$0 < length(I) < \pi$ such that $h_3(\vec{x}) = 0$ if $\vec{x} \in I$ and
$h_3(\vec{x}) = 1$ if $\vec{x} \notin I$, that is, $S_f$ is `$\SS^2$ minus two
spherical triangles'.

As $S_f$ is invariant under every $\varphi_c^{ort}$, vertices of
the triangles are mapped to vertices. Hence we may assume (by splitting $C$
into six pieces) that $(0,0,1)$ is fixed. But then the triangles are also
fixed sets, and every $\varphi_c^{ort}$ is the identity, so we are done as
in the previous case.

This finishes the proof of the Main Theorem.
\ep

\section{Open questions}
\lab{s:open}

\bq
\lab{q:meas}
In the Main Theorem can we relax the assumption of continuity to Lebesgue
measurability, Baire measurability, Borel measurability, Baire class one,
separate continuity or at least one point of continuity?
\eq

\bq
 Which notion of largeness of $C$ suffices for the various results of this
 paper? For example, does the Main Theorem hold if we only assume that $C$
 contains three elements that pairwise generate dense multiplicative subgroups
 of $(0,\infty)$?
\eq

\br
It was shown in
 \cite{Ri} that two such elements suffice for the analogous one-variable
 result. However, two independent elements are not enough here, since if $g$
 is vertically rigid for $c_1$ via a translation and $h$ is vertically rigid
 for $c_2$ via a translation then
 $f(x,y) = g(x) h(y)$ is vertically rigid for both.

 Moreover, the main point in that proof in \cite{Ri} is to replace `splitting
 $C$' by alternative arguments, and we were unable to do so here.
\er

The following question is rather vague.

\bq
\lab{q:rig}
Let us call a set $H \su \SS^2$ rigid if $\psi_c(H)$ is isometric to $H$ for
every $c > 0$. Is there a simple description of rigid sets? Or if we assume
some regularity?
\eq

And finally, the most intriguing problem.

\bq
What can we say if there are more than two variables?
\eq





\end{document}